\documentclass[12pt,oneside, reqno]{amsart}
 \usepackage{graphicx}
 \usepackage{hhline}
\usepackage{amsmath,amsthm,amsfonts,amscd,amssymb,mathrsfs}
\usepackage{xspace}
\usepackage[all]{xypic}
\usepackage{booktabs} % for nicer tables
\usepackage{physics}

\usepackage{array}   % for \newcolumntype macro
\newcolumntype{C}{>{$}c<{$}} % math-mode version of "c" column type
\newcolumntype{L}{>{$}l<{$}} % math-mode version of "l" column type
\newcolumntype{R}{>{$}r<{$}} % math-mode version of "r" column type

\usepackage[backref]{hyperref}
\usepackage[backrefs,initials]{amsrefs}
\usepackage{verbatim}
\usepackage{amscd}   % for commutative diagrams
\usepackage[all]{xy} % for complicated commutative diagrams
\usepackage{youngtab} % for Young tableaux
\usepackage{ytableau} %for Young Tableaux

\usepackage{nicefrac}
\usepackage{xfrac}

\usepackage{mathdots}
\usepackage{tikz}

\usepackage[T1]{fontenc} %% for spacing of << and >>
\usepackage{cleveref} %for multiple references in one and the command \cref
\numberwithin{equation}{section} % number equations within section
\usepackage{color}

\begin{document}

\newcommand{\ccircle}[1]{* + <1ex>[o][F-]{#1}}
\newcommand{\ccirc}[1]{\xymatrix@1{* + <1ex>[o][F-]{#1}}}
%% for re-using a theorem number
\makeatletter
\newtheorem{rep@theorem}{\rep@title}
\newcommand{\newreptheorem}[2]{%
\newenvironment{rep#1}[1]{%
 \def\rep@title{#2 \ref{##1}}%
 \begin{rep@theorem}}%
 {\end{rep@theorem}}}
\makeatother
\newtheorem{theorem}{Theorem}[section]
\newreptheorem{theorem}{Theorem}
\newreptheorem{lemma}{Lemma}
\newtheorem{conjecture}[theorem]{Conjecture}
\newtheorem{lemma}[theorem]{Lemma}
\newtheorem{lem}[theorem]{Lemma}
\newtheorem{prop}[theorem]{Proposition}
\newtheorem{obs}[theorem]{Observation}
\newtheorem{notation}[theorem]{Notation}
\newtheorem{cor}[theorem]{Corollary}
\newtheorem{question}[theorem]{Question}
\theoremstyle{definition}
\newtheorem{definition}[theorem]{Definition}
\newtheorem{examplex}{Example}
\newenvironment{example}
  {\pushQED{\qed}\renewcommand{\qedsymbol}{$\diamondsuit$}\examplex}
  {\popQED\endexamplex}
\newtheorem{xca}[theorem]{Exercise}
\theoremstyle{remark}
\newtheorem{remark}[theorem]{Remark}

\newcommand{\defi}[1]{\textsf{#1}} % for defined terms
\newcommand{\isom}{\cong}
\newcommand{\im}{\operatorname{im}}
\newcommand{\Id}{\text{Id}}
\newcommand{\pr}{\text{pr}}
\newcommand{\Proj}{\operatorname{Proj}}
\newcommand{\Hom}{\operatorname{Hom}}
\newcommand{\End}{\operatorname{End}}
\newcommand{\Gr}{\operatorname{Gr}}
\newcommand{\SL}{\operatorname{SL}}
\newcommand{\GL}{\operatorname{GL}}
\newcommand{\GG}{\operatorname{\text{GME}^{grass}}}
\newcommand{\A}{\mathcal{A}}
\newcommand{\B}{\mathcal{B}}
\newcommand{\C}{\mathcal{C}}
\newcommand{\K}{\mathcal{K}}
\newcommand{\E}{\mathcal{E}}
\newcommand{\F}{\mathcal{F}}
\newcommand{\R}{\mathcal{R}}
\newcommand{\J}{\mathcal{J}}
\newcommand{\I}{\mathcal{I}}
\newcommand{\V}{{\mathcal V}}

\newcommand{\PP}{\mathbb{P}}
\newcommand{\FF}{\mathbb{F}}
\renewcommand{\SS}{\mathbb{S}}
\newcommand{\TT}{\mathbb{T}}
\newcommand{\RR}{\mathbb{R}}
\newcommand{\NN}{\mathbb{N}}
\newcommand{\CC}{\mathbb{C}}
\newcommand{\ZZ}{\mathbb{Z}}
\newcommand{\QQ}{\mathbb{Q}}
\newcommand{\Sym}{\operatorname{Sym}}
\newcommand{\Spec}{\operatorname{Spec}}
\newcommand{\Chow}{\operatorname{Chow}}
\newcommand{\Seg}{\operatorname{Seg}}
\newcommand{\Sub}{\operatorname{Sub}}
\def\bw#1{{\textstyle\bigwedge^{\hspace{-.2em}#1}}}
\def\o{ \otimes }
\def\phi{ \varphi }
\def\ep{ \varepsilon}
\def \P{\mathcal{P}}
\def \U{\mathcal{U}}
\def \S{\mathfrak{S}}
\def \F{\mathcal{F}}
\def \a{\alpha}
\def \b{\beta}
\def \n{\mathfrak{n}}
\def \h{\mathfrak{h}}
\def \d{\mathfrak{d}}
\def \z{\mathfrak{z}}
\def \fb{\mathfrak{b}}
\def \c{\mathfrak{c}}
\def \s{\mathfrak{s}}
\def \ga{\gamma}

\def \g{\mathfrak{g}}
\def \fa{\mathfrak{a}}
\def \e{\mathfrak{e}}
\def \gl{\mathfrak{gl}}
\def \sl{\mathfrak{sl}}
\def \sp{\mathfrak{sp}}
\newcommand{\Ad}{\operatorname{Ad}}
\newcommand{\ad}{\operatorname{ad }}

\newcommand{\SYT}{\operatorname{SYT}}
\newcommand{\Mat}{\operatorname{Mat}}
\newcommand{\SSYT}{\operatorname{SSYT}}
\renewcommand{\span}{\operatorname{span}}
\newcommand{\codim}{\operatorname{codim}}
\newcommand{\diag}{\operatorname{diag}}

\newcommand{\brank}{\operatorname{BRank}}
\newcommand{\Frank}{\operatorname{F-Rank}}
\newcommand{\Prank}{\operatorname{P-Rank}}

\newcommand{\ch}[7]{#1\,#2\,#3\,#4\,#5\,#6\,#7}

\newcounter{nameOfYourChoice}

\def\red#1{{\textcolor{red}{#1}}}
\def\blue#1{{\textcolor{blue}{#1}}}
\def\white#1{{\textcolor{white}{ #1}}}
\newcommand{\arxiv}[1]{\href{http://arxiv.org/abs/#1}{{\tt arXiv:#1}}}

\date{\today}

\author{Luke Oeding}\email{oeding@auburn.edu}
\address{Department of Mathematics and Statistics,
Auburn University,
Auburn, AL, USA
}

\title[\texttt{E\MakeLowercase{xterior}E\MakeLowercase{xtensions}}: A package for \texttt{M2}]{\texttt{E\MakeLowercase{xterior}E\MakeLowercase{xtensions}}: A package for Macaulay2}
\begin{abstract} 
We explain a Macaulay2 implementation of a construction, which appeared in \cite{HolweckOeding23}, of a graded algebra structure on the direct sum of a Lie algebra $\mathfrak{g}$ (typically $\sl_n$) and a $\mathfrak{g}$-module (typically a subspace of an exterior algebra $\bw \bullet\mathbb{C}^n$). We implement brackets, a Killing form, matrix representations of adjoint operators, and ranks of blocks of their powers.
\end{abstract}
\maketitle

\section{Introduction}

The Lie algebra $\g = \sl_n(\FF)$ is a matrix algebra over a field $\FF$ with bracket $[A,B] = AB-BA$ for traceless matrices $A,B \in \FF^{n\times n}$. The package \texttt{ExteriorExtensions} extends this algebra to a larger one $\fa = \g \oplus M$ for a $\g$-module $M$. At this time we have only implemented this for $M$ a sum of exterior powers of the underlying vector space $V = \FF^n$. We say that the extension algebra has GJD (Jordan decomposition consistent with the $G$-action) if the $G = \text{Lie}(\g)$-action on $M$ is consistent with conjugation on adjoint operators $\Ad_T\colon \fa \to \fa$ (see \cite{HolweckOeding23}) for $T\in \fa$. Our development is focused on making tools for studying elements $x \in M$. In particular, the conjugation invariants of an adjoint operator $\Ad_x$, such as the rank, the trace, and characteristic polynomial, are also $G$ invariants when the algebra $\fa$ has GJD. In addition, the adjoint operators inherit a block structure from the grading on $\fa$. 
This package enables one to compute these invariants for many situations of interest in particular for the invariant theory of tensors, and may be useful for applications in quantum information and algebraic statistics, as well as for the study of geometric invariants of special varieties like hyper K\"ahler manifolds, K3 surfaces and more. We explained these connections in \cite{HolweckOeding23}.

The main functions of the package are the following: \texttt{buildAlgebra}, which constructs a new type of graded ring, called an \texttt{ExteriorExtension}, that includes products implemented as \texttt{bracket}, and adjoint operators \texttt{ad}. These rings represent the algebra $\fa$. The package computes standard matrix representations of adjoint operators on $\fa$, and facilitates the tabulation of the ranks of blocks of powers of adjoint operators via \texttt{blockRanks}. We also compute the \texttt{KillingMatrix}, which may help identify the algebra $\fa$. 

The remainder of the introduction provides a brief mathematical account of each of these objects. In Section~\ref{sec:examples} we provide examples of the functionality of the package, specifically these functions, and comment on the meaning of the computations.  We note that \texttt{M2} starts numbering by 0 by default, whereas we often start numbering at 1, and we hope that this shift between explanations and implementations is not too confusing for the reader.

\subsection{Exterior algebras} The \defi{exterior algebra} on a vector space $V$ is the direct sum of exterior powers $\bw\bullet V = \bigoplus_{k\geq 0} \bw kV$, with product given by concatenation, which respects a $\ZZ$-grading. The exterior algebra on $V$ is Artinian when $V$ is finite-dimensional.  However, we can obtain a $\ZZ_n$-grading if $n = \dim V$ by choosing a volume form $0\neq \Omega \in \bw n V$ and utilizing the correspondingly defined Hodge star $\star$ operation to replace products that would have degree greater than $n$. We explain this more fully in Section~\ref{sec:structureTensor}. 

The vector space $V$ carries the standard group action by $G = \SL(V)$ and Lie algebra action by $\mathfrak{g} = \sl(V)$. These actions are induced up to the exterior powers, making any direct sum of exterior powers $G$ (or $\mathfrak g$)-modules.

\subsection{Extending the Lie algebra by appending an exterior algebra} 
Now fix positive integers $k$ and $n = \dim V$. The exterior algebra $\bw\bullet V = \bigoplus_{k\geq 0} \bw kV$ is graded by degree, however one obtains a (possibly) different grading by declaring $\bw j V$ to be the grade-1 piece, and defining grading on the rest of the graded pieces inductively.
More specifically, consider the graded vector space 
\[\fa = \g \oplus \bw{k} V \oplus \bw{2k} V \oplus \dots  = \g \oplus \bigoplus_{ik\not \equiv 0 \mod n} (\bw {ik} V),
\]
where we replace the scalars from the exterior algebra with $\fa_0 = \g$ and $\fa_i = \bw{ik \mod n} V$ for each $i$ from $i=1$ up to but not including the first $i$ such that  $ik\equiv 0 \mod n$. Note that if $k$ divides $n$ we are appending part of a sub-algebra of the exterior algebra of the parts whose degrees are multiples of $k$. For example, when $k=3$, and $n=9$ we construct 
\[\fa = \sl_9 \oplus \bw 3 \FF^9 \oplus \bw 6 \FF^9,\]  
but  when $k=3, n=10$ we construct
\[\fa = \sl_{10} \oplus \bw 3 \FF^{10} \oplus \bw 6 \FF^{10}\oplus \bw 9 \FF^{10} \oplus \bw 2 \FF^{10}\oplus \bw 5 \FF^{10}\oplus \bw 8 \FF^{10}\oplus \bw 1 \FF^{10}\oplus \bw 4 \FF^{10}\oplus \bw 7 \FF^{10}.\] 

This vector space $\fa$ is constructed by the function \texttt{buildAlgebra} and the pieces that comprise it are accessed by \texttt{LieAlgebra} (for $\fa_0$) and \texttt{appendage} (for $\fa_{>0}$). More specifically, we construct $\fa_0 = \sl(V)$ by making a common choice basis for the traceless $n \times n$ matrices, and we access $\fa_{>0}$ by calling a ring with skew-commuting variables, $e_i$. One can access the bases we choose for these algebras by \texttt{bases}, which produces a hash table with key \texttt{i} calling the basis of $\fa_i$.

Already the graded vector space  space $\fa$ has brakets (products) defined for $\fa_i \times \fa_j \to \fa_{i+j}$ when $1<i,j,i+j<n$ (the skew-commuting product on the ring $\FF[\FF^n]$). In the next section, we explain how to extend the definition of the bracket to the entire graded vector space in an equivariant fashion. 

One of the main results from \cite{HolweckOeding23} tells precisely when it is possible to make $\fa$ into a graded algebra respecting the group structure.
\begin{theorem}{\cite[Theorem~3.18]{HolweckOeding23} }
The vector space $\fa = \sl_{n} \oplus \bigoplus_{k=1 \ldots n-1}\bw{k} \CC^{n}$ has a $\ZZ_n$-graded algebra structure with a Jordan decomposition consistent with the $G = \SL(V)$-action. There is a unique (up to scale) equivariant bracket product that agrees with the $\g$-action on each $\g_i$. If $n = 2k$, the equivariant bracket must satisfy the property that the restriction to $\bw{k} \CC^{n} \times \bw{k} \CC^{n} \to \g$ must be commuting when $k$ is odd and skew-commuting when $k$ is even. For any $k$ such that $2k \neq n$ the bracket $\bw{k} \CC^{n} \times \bw{k} \CC^{n} \to \bw{2k \mod n}\CC^n$ must be skew-commuting when $k$ is odd and commuting when $k$ is even. 
\end{theorem}
This result here is that a natural equivariant bracket exists, so it makes sense to try to implement it in the first place, and moreover, since the structure tensor $B$ for the algebra $\fa$ is essentially unique, there aren't many choices in its formulation. In particular, the only choices we have are choices of scalar multiples of the blocks $B_{i,j,i+j}$ of the tensor. 
We implement that product as a method function in \texttt{M2} called \texttt{bracket}. Next, we give a detailed description of the essentially unique structure tensor for this algebra.

\subsubsection{The structure tensor defining the bracket}\label{sec:structureTensor}
 %Here are some additional details on our implementation of the bracket product for an exterior extension.
If $V$ is a vector space over a field $\FF$, we let $V^*$ denote the dual vector space of $\FF$-linear functionals. Following the perspective of \cite{LandsbergTensorBook} a tensor $T \in V_1\otimes \cdots \otimes V_m$ can be equivalently viewed in many ways, such as a linear functional $T\colon V_1^*\otimes \cdots \otimes V_m^* \to \FF$, as a linear map $T\colon V_1^*\otimes \cdots \otimes \widehat{V^*_i} \otimes\cdots \otimes V_m^* \to V_i$ (where $\widehat {\cdot}$ denotes omission), or as a multilinear map $T\colon V_1^*\times \cdots \times V_m^* \to \FF$.

Recall that a \defi{structure tensor} $B \in \fa^* \otimes \fa^* \otimes \fa$ representing the bilinear product $[\;,\;] \colon \fa \times \fa \to \fa$  is equivalent information to the algebra structure on $\fa$. We denote by $B_{i,j,k} \in  \fa_i^* \otimes \fa_j^* \otimes \fa_k$ the graded pieces of $B$. To respect the grading we must require that $B_{i,j,k} = 0$ whenever $i+j \not\equiv k \mod n$. 

For $B_{0,0,0}$ we have the standard skew-commuting bracket on $\sl_n$, namely $[A,B] = AB-BA$. 
For $B_{0,j,j}$ (and for $ B_{j,0,j}$) with  $j>0$ we use the standard Lie algebra action of $\sl_n$ on the fundamental representation $\bw{jk}\CC^n$, namely for $A \in \sl_n$ and $T \in \bw{jk}V$ we take $[A,T] := A.T$, where $.$ denotes the Lie algebra action. Similarly $[T,A] = -A.T$.

For $i>0$  we have two possibilities for $B_{i,j,(i+j) \mod n}$  depending on whether or not $(i+j) \mod n=0$. When $(i+j) \equiv 0 \mod n $ we separate two further cases depending on whether $i\leq n/2$.  
When $0<i+j <n$ we use the standard product in the exterior algebra as the bracket. 

If $i+j \equiv 0 \mod n$, for $S$ of degree $i\leq n/2$ and $T$ of degree $j> n/2$ we utilize the Hodge star $\star T$ to change degrees and compute as follows. We declare an orthonormal basis $e_1,\ldots,e_n$ of $V = \FF^n$ and note that the Hodge star is determined by a volume form which we choose to be $\Omega = e_1\wedge\cdots \wedge e_n$.
The standard Hodge star operator for $T \in \bw j V $  producing $\star T \in \bw{n-j} V$ is implemented as follows:
\[
\star T =
\sum_{\alpha \in \binom{[n]}{j}} \text{sgn}(\alpha)\left(\big(\prod_{i \in \alpha} e_i\big)\lrcorner T  \right )\prod_{i\in \alpha^c} e_i,
\]
where $\binom{[n]}{j}$ denotes the ordered subsets of size $j$ of $[n]=\{1,\ldots,n\}$, the products are taken in the exterior algebra $\bw \bullet [e_1,\ldots,e_n]$, and $\lrcorner$ denotes contraction. We note that while it is certainly the case that $\bw{n-j}V \cong \bw{j}V^*$ (the dual vector space), for $T \in \bw j V$ we still consider $\star T \in \bw{n-j}V$ in the not-dual coordinates since in this case $i+j =  n$ with $n>j>n/2$ we have $n-j = i$, so the pair $(S,\star T)$ consists of elements of equal degree ($i$).
Taking a pair of partial derivatives with respect to elements of degree 1 takes a pair $(S,\star T)$ of elements of degree $i$ to a pair $\left( \frac{\partial S}{\partial e_p} , \frac{\partial \star T}{\partial e_q} \right ) $ of elements both of degree $i-1$. Contracting the resulting pair produces a number, which we store as the $(p,q)$-entry of a matrix. Projection (denoted $\pi$) to traceless matrices produces an element of $\sl_n$,
\[\begin{matrix}
\pi\colon \; \FF^{n\times n} &\to& \sl_n \\
A &\mapsto& A- \frac{\tr(A)}{n}I_{n\times n} \;.
\end{matrix}
\]
We have chosen a global scalar to calibrate and match the outcomes of \cite{Vinberg-Elashvili} (see the end of example~\ref{ex:e8}) which considered the case $\sl_9 \oplus \bw{3}\CC^9\oplus \bw{3}\CC^9$ as a model for the exceptional Lie algebra $\mathfrak{e}_8$.
 In summary, the map that takes two elements of the exterior extension to a matrix in $\sl_n$ is the following:
\[ 
\begin{matrix}
\fa_i \times \fa_j &\to& \fa_0 \\
(S,T) & \mapsto & [S,T] =  -k \pi \left(  \left( \frac{\partial S}{\partial e_p} \lrcorner \frac{\partial \star T}{\partial e_q} \right )_{1\leq p,q\leq n} \right ) 
\end{matrix} \quad \quad \text{when } i+j = n \text{ and } j>n/2.
\]

For $i\geq n/2$ and $j< n/2$ we reverse the role of the Hodge dual and proceed as in the previous paragraph.

Finally, when $i+j>n$ we use Hodge star to convert to a product that we have already computed. If $S$ has degree $i>n/2$ or $T$ has degree $j>n/2$ then we take $[S,T] = \star [\star S,\star T]$ and use the bracket previously defined for $i+j\leq n$. 

Finally, we remark that the products that have degree sum $i+j>n$ on the level of basic elements $S = e_I$ and $T = e_J$ essentially result in $[S,T] \sim e_{I\cap J}$, however one has to pay attention to signs. Likewise, when the degree sum is $i+j<n$ on the level of basic elements $S = e_I$ and $T = e_J$ essentially results in $[S,T] \sim e_{I\cup J}$, however again one has to pay attention to signs.

\subsubsection{Adjoint Operators}
The standard matrix representation of the adjoint operator on $\fa$ defined by $\Ad_x(y)  = [x,y]$ is constructed by choosing an ordered basis of $\fa$, and expressing the results of $\Ad_x(y)$ for each basis vector $y$ in terms of the standard basis as the corresponding column.
One caveat is that we have to choose a basis of $\sl_n$, and there is not a "standard" choice for the basis of the Cartan subalgebra of traceless diagonal matrices. We choose a common one $\{ E_{ii}-E_{i+1,i+1} \mid i=1\ldots n-1\}$. Then we have to solve a system of equations to find the coordinates for the cartan rather than just extracting the entries of the matrix. Solving this system of equations can depend on the field $\FF$ because we invert a matrix over $\FF$ to find the solution.

Not that the adjoint operator is a linear map $\Ad_x \colon \fa \to \fa$, which is represented by a (large) matrix. The structure tensor is graded so that $B_{i,j,k} \colon \fa_i \times \fa_j \to \fa_k$. The map respects the grading on $\fa$ implies that $B_{i,j,k} = 0 $ unless $k \equiv i+j \mod n$.   This, in turn, induces a block structure on the matrix representing $\Ad_x$ since $\Ad_x(y) = \bigoplus_{i,j,k} (B_{i,j,k}(x,y))$.
 For example, in the 2-graded case, with $x \in \fa_1$ we have the blocking
\[
\Ad_x = \left( \begin{smallmatrix} 0 & * \\ * & 0 \end{smallmatrix} \right),
\]
where the $0$'s and $*$'s respectively represent block matrices, and column blocks (and row blocks) are respectively of sizes $(\dim \fa_0$, $\dim \fa_1)$. 
In the 3-graded case, with $x \in \fa_1$ we have the blocking
\[
\Ad_x = \left( \begin{smallmatrix}  0 & 0 & * \\ * & 0 & 0\\  0 & * & 0   \end{smallmatrix} \right),
\]
with column blocks (and row blocks) respectively of sizes $(\dim \fa_0$, $\dim \fa_1, \dim \fa_2)$. 

 The ranks of the blocks of $\Ad_x$ are invariants of $x$, as are the ranks of blocks of the powers of $\Ad_x$. We collect these ranks in a table and display them with a function \texttt{prettyBlockRanks}. 
 
 We note that we implement the ranks of powers of blocks ``under the hood'' by extracting blocks of $(\Ad_x)^k$ as submatrices of the appropriate sizes (determined by the dimensions of the graded pieces, and the respective ordered bases we constructed to carry out our computations) and computing their ranks, and moreover, we compute $(\Ad_x)^k$ sequentially.  We note that it is, in theory, possible to make a generic adjoint operator since the map $\fa \to \End \fa$ is linear, and that would produce $\Ad_x$ as a matrix of linear forms on $\fa$. Computing with a large matrix of linear forms can quickly become costly computationally, so we prefer to only work with specific numerical instances, where computing powers and ranks of powers is computationally inexpensive.

The Jacobi identity is stated via adjoint operators as $[\Ad_x, \Ad_y] = \Ad_{[x,y]}$ \cite{FultonHarris}. We can check whether this identity holds for the algebras we construct either on random elements of each grade, or on bases.

\subsubsection{The Killing form}  Recall that the Killing form for an algebra is a quadratic form on $\fa$ defined by $B(x,y) = \trace (\Ad(x)\cdot \Ad(y))$ (unfortunately the symbol $B$ is being re-used.)  A \defi{Lie algebra} is an algebra that has a skew-commuting bracket that satisfies the Jacobi identity. Cartan's criterion states that a Lie algebra is semi-simple (a direct sum of simples) if and only if the Killing form is non-degenerate. Angles between weight basis vectors are obtained via the inner product induced from the Killing form, which in turn defines a Dynkin diagram. This information can be used to identify the Lie algebra.

\subsubsection{On fields of definition}This package was written with the fields $\CC$ or $\RR$ in mind because of our knowledge of results in representation theory over the real and complex numbers. However, much of the functionality of \texttt{M2} is better suited for exact arithmetic, over $\QQ$ or $\ZZ_p$, for instance. In one place the package attempts to find an expression of a matrix in terms of a specially chosen basis of $\sl_n$, which we carry out by solving a system of equations. Solving this linear system inverts a matrix (the locally defined matrix \texttt{Qmat} in the package), and in \texttt{M2} this operation has different behavior based on the ground field. At this time of development, the author has not tested behavior over other rings or fields, and so the functions are implemented over the \texttt{M2} type \texttt{Ring}, and \texttt{M2} will throw an error when imprecise fields like $\RR$ or $\CC$ are passed to the \texttt{buildAlgebra} command. 
It is still possible to have some of this functionality, and to do this is typically best to start working over an exact field, and only lift to $\RR$ or $\CC$ when necessary or when it makes sense to do so. By not working over inexact fields until necessary we reduce the chances for numerical errors.

\section{Examples}\label{sec:examples}
The algebras $\fa = \sl_{2m}\oplus \bw m \CC^{2m}$ for $m=2$ and $m=4$ are actually simple Lie algebras, respectively $\sl_4 \oplus \bw{2}\CC^4 \cong \sp_6$ (a symplectic Lie algebra) and $\sl_8 \oplus \bw{4}\CC^8 \cong \mathfrak{e}_7$ (an exceptional Lie algebra).  An excellent overview of fine gradings of simple Lie algebras is given in \cite{MR3601079}. It is not immediately clear to the author if this example falls within those fine gradings classified in \cite[Example~6]{MR3601079}.

\begin{example} In this example we show the use of the functions following, \texttt{buildAlgebra}, \texttt{makeTraceless}, \texttt{bracket}, \texttt{ad}, \texttt{appendage}, and \texttt{KillingMatrix}, to illustrate the known sporadic isomorphism  $\sl_4 \oplus \bw{2}\CC^4 \isom \sp_6$. 
Recall from \cite[Ch. 15.2 and Ch.17.1]{FultonHarris} the structure of the two algebras. The weights of $\bw 2\CC^4$ are the vertices of an octahedron, while the weights of $\sl_4$ are the midpoints of the edges of a cube. These geometries are also present in the weights of a basis of $\sp_6$, see \cite[p220, p254]{FultonHarris}. 

In this example we construct $\fa = \sl_4(\FF) \oplus \bw2 \FF^4$ for  $\FF =\QQ$, and note some of the salient features of this algebra extension, even though we may be mainly interested in the theory of the complex Lie algebras. Near the end of the example, we note that the Killing form of the algebra has full rank and all real eigenvalues, and this means that the isomorphism $\fa \isom \sp_6$  holds as in isomorphism of real Lie algebras.

We set the algebra as an \texttt{ExteriorExtension} via:
\begin{verbatim}
     loadPackage"ExteriorExtensions"
     ea24 = buildAlgebra(2,4,e,QQ)
\end{verbatim}

We will see that $\sl_4 \oplus \bw 2 \CC^4$ with the (essentially unique) bracket coming from the GJD criterion forms a Lie algebra. We can check (at least on random points with rational coefficients) if the bracket satisfies the Lie algebra axioms in each graded piece of the algebra. If we had the patience we could also check these identities for all pairs of basis vectors. For illustration we only perform tests for randomly chosen elements.  For grade 0 times grade 0, we first make 2 random traceless matrices:
\begin{verbatim}
     A0 = makeTraceless random(QQ^4,QQ^4)
     B0 = makeTraceless random(QQ^4,QQ^4)
\end{verbatim}

To access the \texttt{bracket} and adjoint operator \texttt{ad} for the algebra \texttt{ea24} that we have constructed, we do so with the dot, which we see in the following random check of the skew symmetry and Jacobi identities:
\begin{verbatim}
     ea24.bracket(A0,B0) + ea24.bracket(B0,A0)
          0
     ea24.ad(ea24.bracket(A0,B0)) - ea24.bracket(ea24.ad(A0),ea24.ad(B0))
          0
\end{verbatim}

\noindent We shouldn't be surprised that we get 0 for both since the grade 0 piece is a Lie algebra.

For a pair of grade 1 elements we call for random elements of the \texttt{appendage} and check their brackets are skew and satisfy the Jacobi identity as follows:
\begin{verbatim}
     A1 = random(2, ea24.appendage)
     B1 = random(2, ea24.appendage)
     ea24.bracket(A1,B1) + ea24.bracket(B1,A1)
          0
     ea24.ad(ea24.bracket(A1,B1)) - ea24.bracket(ea24.ad(A1),ea24.ad(B1))
          0
\end{verbatim}

    For a point in grade 0 and a point in grade 1, we see the axioms hold, but we recall the caveat that it is important to make the matrix traceless (we implemented this projection as \texttt{makeTraceless}), otherwise the Jacobi identity fails. We proceed as before:
\begin{verbatim}
     ea24.bracket(A0,B1) + ea24.bracket(B1,A0)
          0
     ea24.ad(ea24.bracket(A0,B1)) - ea24.bracket(ea24.ad(A0),ea24.ad(B1))
          0
\end{verbatim}
We can also construct a Killing matrix and see that it is \(21 \times 21\) and non-degenerate. The only possibility for a simple Lie algebra over $\CC$ of dimension 21  is \(\mathfrak{sp}_6\).
\begin{small}
    \begin{verbatim}
     K = ea24.KillingMatrix()
          | 20  -10 0   0  0  0  0  0  0  0  0  0  0  0  0  0   0  0   0   0  0   |
          | -10 20  -10 0  0  0  0  0  0  0  0  0  0  0  0  0   0  0   0   0  0   |
          | 0   -10 20  0  0  0  0  0  0  0  0  0  0  0  0  0   0  0   0   0  0   |
          | 0   0   0   0  0  0  0  0  0  10 0  0  0  0  0  0   0  0   0   0  0   |
          | 0   0   0   0  0  0  0  0  0  0  10 0  0  0  0  0   0  0   0   0  0   |
          | 0   0   0   0  0  0  0  0  0  0  0  10 0  0  0  0   0  0   0   0  0   |
          | 0   0   0   0  0  0  0  0  0  0  0  0  10 0  0  0   0  0   0   0  0   |
          | 0   0   0   0  0  0  0  0  0  0  0  0  0  10 0  0   0  0   0   0  0   |
          | 0   0   0   0  0  0  0  0  0  0  0  0  0  0  10 0   0  0   0   0  0   |
          | 0   0   0   10 0  0  0  0  0  0  0  0  0  0  0  0   0  0   0   0  0   |
          | 0   0   0   0  10 0  0  0  0  0  0  0  0  0  0  0   0  0   0   0  0   |
          | 0   0   0   0  0  10 0  0  0  0  0  0  0  0  0  0   0  0   0   0  0   |
          | 0   0   0   0  0  0  10 0  0  0  0  0  0  0  0  0   0  0   0   0  0   |
          | 0   0   0   0  0  0  0  10 0  0  0  0  0  0  0  0   0  0   0   0  0   |
          | 0   0   0   0  0  0  0  0  10 0  0  0  0  0  0  0   0  0   0   0  0   |
          | 0   0   0   0  0  0  0  0  0  0  0  0  0  0  0  0   0  0   0   0  -10 |
          | 0   0   0   0  0  0  0  0  0  0  0  0  0  0  0  0   0  0   0   10 0   |
          | 0   0   0   0  0  0  0  0  0  0  0  0  0  0  0  0   0  0   -10 0  0   |
          | 0   0   0   0  0  0  0  0  0  0  0  0  0  0  0  0   0  -10 0   0  0   |
          | 0   0   0   0  0  0  0  0  0  0  0  0  0  0  0  0   10 0   0   0  0   |
          | 0   0   0   0  0  0  0  0  0  0  0  0  0  0  0  -10 0  0   0   0  0   |
     rank K
          21
\end{verbatim}
\end{small}
We implemented the \texttt{KillingMatrix} as a function rather than something that gets constructed when the algebra is constructed since, without additional tricks, it may be time-consuming to construct this matrix.
One checks that the Killing matrix has a full set of real eigenvectors, so the isomorphism  \(\mathfrak{a} \isom \mathfrak{sp}_6\) holds over $\RR$.
\begin{verbatim}
     toList eigenvalues K 
          {34.1421, 20, 5.85786, 10, -10, 10, -10, 10, -10, 10, -10, 10,
           -10, 10, -10, 10, -10, 10, -10, 10, -10}
\end{verbatim}
Here we note the first instance of a change in the field of definition. To compute eigenvalues M2 automatically upgraded the matrix $K$, or at least its eigenvalues, to work over $\CC$.

The tensors $e_i\wedge e_j \pm e_k\wedge e_l$ for $i,j,k,l$ distinct have diagonalizable adjoint forms (we call them semi-simple or ad-semi-simple). We can check that the geometric and algebraic multiplicities of the eigenvalues add up:
\begin{verbatim}
     E = ea24.appendage
     A = ea24.ad(E_0*E_1 - E_2*E_3)
     tally apply(toList eigenvalues (A), xx-> round(2, realPart xx) +
     round(2, imaginaryPart xx)*ii)
          Tally{-1.41 => 5}
                0 => 11
                1.41 => 5
     rank A
          10
\end{verbatim}
We also note that even though $A$ is defined over \texttt{QQ}, we can let M2 work with square roots.
\begin{verbatim}
     tally apply(toList eigenvalues (1/sqrt(2)*A), xx->
     round(2, realPart xx) + round(2, imaginaryPart xx)*ii)
          Tally{-1 => 5}
                0 => 11
                1 => 5
\end{verbatim}
The semi-simple structure is also indicated by the block ranks computation, via \texttt{blockRanks} (to produce a matrix as a result) or via \texttt{prettyBlockRanks} for onscreen printing:
\begin{verbatim}
     ea24.prettyBlockRanks A
           +---+---+---+---+-----+
           |g00|g01|g10|g11|total|
           +---+---+---+---+-----+
           |0  |5  |5  |0  |10   |
           +---+---+---+---+-----+
           |5  |0  |0  |5  |10   |
           +---+---+---+---+-----+
           |0  |5  |5  |0  |10   |
           +---+---+---+---+-----+
           |5  |0  |0  |5  |10   |
           +---+---+---+---+-----+
\end{verbatim}
The rows in the above table correspond to subsequent powers of the matrix, and we stop computation after sufficient repetition or zero total rank occurs.

We can check that neither the elements nor their adjoint forms commute, so they're not simultaneously diagonalizable.
\begin{verbatim}
     ea24.bracket(E_0*E_1 + E_2*E_3,E_0*E_2 + E_1*E_3)
          {1} | 0  0 0 -1 |
          {1} | 0  0 1 0  |
          {1} | 0  1 0 0  |
          {1} | -1 0 0 0  |
\end{verbatim}

We checked that the only elements of $\bw 2\CC^4$ that commute with $e_i\wedge e_j + e_k\wedge e_l$ are scalar multiples of itself. So we can't hope for a Cartan subalgebra that intersects $\bw2 \CC^4$ and includes more than one of these forms. We did this computation by accessing the basis (of the grade 1 part) of the algebra via \texttt{ea24.bases}.
\begin{verbatim}
     L = apply(ea24.bases#1, xx-> ea24.bracket(xx,e_0*e_1 + e_2*e_3))
\end{verbatim}
Then we prepare additional variables, and solve the system of equations coming from the commutation conditions:
\begin{verbatim}
     R = QQ[x_0.. x_5]
     decompose ideal flatten entries sum(6, i-> sub(L_i, R)*x_i)
          ideal (x , x , x , x , x  - x )}
                  4   3   2   1   0    5
\end{verbatim}
Here we note that we have changed rings again only when necessary. Attention should be paid to the ring of definition for the matrices the package produces, and sometimes substitution commands are necessary.
\end{example}
 \begin{example}
We consider $\sl_6 \oplus \bw3 \CC^6$, which is not a Lie algebra as it fails to be skew-commuting (it is commuting on the grade 1 piece) and it fails the Jacobi identity. 
\begin{verbatim}
     ea36 = buildAlgebra(3,6,QQ);
     A = random(3, ea36.appendage);
     B = random(3, ea36.appendage);
     ea36.bracket(A,B) + ea36.bracket(B,A) == 0
          false
     ea36.bracket(A,B) - ea36.bracket(B,A) == 0
          true
     ea36.ad(ea36.bracket(A,B)) - ea36.bracket(ea36.ad(A),ea36.ad(B)) == 0
          false
\end{verbatim}
Interestingly, the trace of the 4th power of the adjoint operator still produces the single (degree 4) generator of the invariant ring $\CC[\bw 3 \CC^6]^{\SL_6} \cong \CC[f_4]$ in this case. Computing values of this invariant is accomplished as follows:
\begin{verbatim}
     E = ea36.appendage
     f4 = A-> trace(A^4)
     N = E_0*E_1*E_2 + E_3*E_4*E_5
     f4 ea36.ad(N)
          36     
\end{verbatim}
\vspace{-1.5em}
\end{example} 
 \begin{example}
 Consider the case of $\sl_8 \oplus \bw 4 \CC^8 \cong \mathfrak{e}_7$ considered in \cite{Antonyan, oeding2022}. Antonyan shows that a 7-dimensional Cartan subalgebra is given by 4-forms of type $e_I + e_J$ with the concatenation $I |J $ a permutation of $[8]$, the 7 pairs of indices $I,J$ are the following
  \[
 \begin{matrix}
\{\{1,2,3,4\},\;\{5,6,7,8\}\}, &&
\{\{1,3,5,7\},\;\{6,8,2,4\}\},&&
\{\{1,5,6,2\},\;\{8,4,7,3\}\},\\
\{\{1,6,8,3\},\;\{4,7,5,2\}\},&&
\{\{1,8,4,5\},\;\{7,2,6,3\}\}, &&
\{\{1,4,7,6\},\;\{2,3,8,5\}\}, \\
\{\{1,7,2,8\},\;\{3,5,4,6\}\}.
 \end{matrix}
\] 
We can check the diagonalizability and commutativity as follows. First, we construct the algebra and our candidate for the basis of the Cartan subalgebra.
 \begin{verbatim}
     ea48 = buildAlgebra(4,8,e,QQ);
     list2E = L-> sum(L, l-> product(l, i-> e_(i-1)));
     aa = {{{1,2,3,4},{5,6,7,8}},{{1,3,5,7},{6,8,2,4}},{{1,5,6,2},{8,4,7,3}},
     {{1,6,8,3},{4,7,5,2}},{{1,8,4,5},{7,2,6,3}},{{1,4,7,6},{2,3,8,5}},
     {{1,7,2,8},{3,5,4,6}}};
\end{verbatim}
Then we produce the adjoint operators for each, and compute respectively their sums of the algebraic and geometric multiplicies of eigenvalues to see that they do, in fact, agree, and hence these elements are diagonalizable.
 \begin{verbatim}
     for i to length(aa)-1 do(
     adList_i = ea48.ad(list2E aa_i);
     print(# for xx in eigenvalues adList_i list if 
     	not round(abs(xx)) ==0 then xx else continue, rank adList_i); )
          (66, 66)
          (66, 66)
          (66, 66)
          (66, 66)
          (66, 66)
          (66, 66)
          (66, 66)
\end{verbatim}
Next, we check that these elements mutually pairwise commute by computing all their commutators. We store these results as a matrix whose entries are the ranks of the commutators, and note that we get the zero matrix:
 \begin{verbatim}
     matrix apply(7, i-> apply(7, j-> rank ea48.bracket(adList_i,adList_j)))
          0
\end{verbatim}
\vspace{-1.5em}
 \end{example}
 
 \begin{example}\label{ex:e8}
 Larger algebras are also possible, such as $\mathfrak{e}_8 \cong \sl_9  \oplus \bw 3 \CC^9 \oplus\bw 6 \CC^9 $ from \cite{Vinberg-Elashvili}, and we have done extensive computations in \cite{HolweckOeding23}. Here are some example computations that check that the bracket satisfies the Lie-bracket conditions (skew-symmetry and Jacobi).
We check bracket conditions for pairs of random elements in each grade:
  \begin{verbatim}
     ea = buildAlgebra(3,9)
     E = ea.appendage
     A0 = makeTraceless random(QQ^9,QQ^9);
     B0 = makeTraceless random(QQ^9,QQ^9);
     A1 = random(3, E);
     B1 = random(3, E);
     A2 = random(6, E);
     B2 = random(6, E);
     ea.bracket(A0,B0) + ea.bracket(B0,A0)
          0
     ea.ad(ea.bracket(A0,B0)) -ea.bracket(ea.ad(A0),ea.ad(B0)) 
          0
     ea.bracket(A0,B1) + ea.bracket(B1,A0)
          0
     ea.ad(ea.bracket(A0,B1)) -ea.bracket(ea.ad(A0),ea.ad(B1)) 
          0
     ea.bracket(A0,B2) + ea.bracket(B2,A0)
          0
     ea.ad(ea.bracket(A0,B2)) -ea.bracket(ea.ad(A0),ea.ad(B2)) 
          0
     ea.bracket(A1,B1) + ea.bracket(B1,A1)
          0
     ea.ad(ea.bracket(A1,B1)) -ea.bracket(ea.ad(A1),ea.ad(B1)) 
          0
     ea.bracket(A1,B2) + ea.bracket(B2,A1)
          0
     ea.ad(ea.bracket(A1,B2)) -ea.bracket(ea.ad(A1),ea.ad(B2)) 
          0
     ea.bracket(A2,B2) + ea.bracket(B2,A2)
          0
     ea.ad(ea.bracket(A2,B2)) -ea.bracket(ea.ad(A2),ea.ad(B2))           
          0
\end{verbatim} 
The following match equations 2.5 and 2.6 from \cite{Vinberg-Elashvili}.
\begin{verbatim}
     pt = E_0*E_1*E_2;
     tmp  = ea.bracket(pt, ea.star(pt)) 
          {2} | -2/3 0    0    0   0   0   0   0   0   |
          {2} | 0    -2/3 0    0   0   0   0   0   0   |
          {2} | 0    0    -2/3 0   0   0   0   0   0   |
          {2} | 0    0    0    1/3 0   0   0   0   0   |
          {2} | 0    0    0    0   1/3 0   0   0   0   |
          {2} | 0    0    0    0   0   1/3 0   0   0   |
          {2} | 0    0    0    0   0   0   1/3 0   0   |
          {2} | 0    0    0    0   0   0   0   1/3 0   |
          {2} | 0    0    0    0   0   0   0   0   1/3 |
     ea.bracket(tmp, pt) 
          -2e e e
             0 1 2
\end{verbatim} 
\vspace{-1.5em}
 \end{example}
\section{The Future}What we anticipate for the future of this package is to provide functionality for theoretical work on tensors. For example, it is possible to compute an adjoint form of a tensor of small format in a short amount of time on a laptop computer. From here one may use the adjoint form to compare tensors via their adjoint block ranks, adjoint spectra, etc. This can have applications whenever one wants to separate tensor orbits, for example. This strategy was explored in \cite{HolweckOeding23}.  The author has work in progress with collaborators where this strategy is applied to tensors that arise in phylogenetic networks. In particular, we are very interested in the cases beyond the ``tame'' tensor formats, where there are not finitely many orbits, and the algebra we construct $\fa$ is not a semi-simple Lie algebra. One very interesting case is $\fa = \sl_{12} \oplus \bw 4 \CC^{12}\oplus \bw 8 \CC^{12}$, which contains a graded subalgebra  \[\fa_{4,4,4} = (\sl_{4} \oplus \sl_{4} \oplus \sl_{4} ) \oplus (\CC^4 \otimes \CC^4 \otimes  \CC^4 ) \oplus (\bw 2 \CC^4 \otimes \bw 2 \CC^4 \otimes \bw 2  \CC^4 )\oplus (\bw 3 \CC^4 \otimes \bw 3 \CC^4 \otimes \bw 3  \CC^4 ),\]
which allows us to utilize these tools to study tensors of format $4\times 4\times 4$. This format is the smallest cubic format for tensors where there is not a complete irredundant classification of orbits. We are quite interested in studying this case further with these new tools.

We have tried to be as efficient as possible in the implementation of the functions in this package, but there is certainly room for improvement, both mathematically and algorithmically. On the mathematical side, if one knows already how to compute the Killing form in a combinatorial fashion it might be possible to skip the lengthy process that computes each entry of the Killing matrix as the trace of a product of adjoint operators.  Algorithmically, it might also be possible to speed up the computation of the adjoint operators in the first place, and it also seems very possible to utilize sparse arrays to be more efficient with memory. In addition, it could be quite interesting to attempt to construct similar extensors for other Lie algebras (replacing $\fa_0 = \sl_n$), or for other modules, replacing $\fa_{<0}$. 

\section{Acknowledgements}
The author thanks Frederic Holweck for his input on the development of this package, which was used for their joint work. Ian Tan also provided useful discussions and testing. The author thanks two anonymous referees for their comments which greatly improved this manuscript. The author acknowledges partial support from the Thomas Jefferson Foundation. 
\bibliography{ /Users/lao0004/Library/CloudStorage/Box-Box/texFiles/main_bibfile.bib}
%\bibliography{/Users/oeding/Library/CloudStorage/Box-Box/texFiles/main_bibfile.bib}

\end{document}